# Stability Analysis of Fractional Order Systems Described in the Lur'e Structure

Shima Sadat Mousavi, *Student Member, IEEE*, and Mohammad Saleh Tavazoei, *Member, IEEE*

Department of Electrical Engineering,

Sharif University of Technology

shimasadat_mousavi@ee.sharif.edu, tavazoei@sharif.edu

*Abstract*—Lur'e systems are feedback interconnection of a linear time-invariant subsystem in the forward path and a memoryless nonlinear one in the feedback path, which have many physical representatives. In this paper, some classical theorems about the $L_2$ input-output stability of integer order Lur'e systems are discussed, and the conditions under which these theorems can be applied in fractional order Lur'e systems with an order between 0 and 1 are investigated. Then, application of circle criterion is compared between Lur'e systems of integer and fractional order using their corresponding Nyquist plots. Furthermore, applying Zames-Falb and generalized Zames-Falb theorems, some classes of stable fractional order Lur'e systems are introduced. Finally, in order to generalize the off-axis circle criterion to fractional order systems, another method is presented to prove one of the theorems used in its overall proof.

## I. Introduction

Many real-world systems can be modeled as closed loop systems which are an interconnection of a linear time-invariant (LTI) subsystem in the forward path with a memoryless nonlinear subsystem in the feedback path (see Figure 1). Such systems are called "Lur'e systems", the name of the first who introduced them [1].

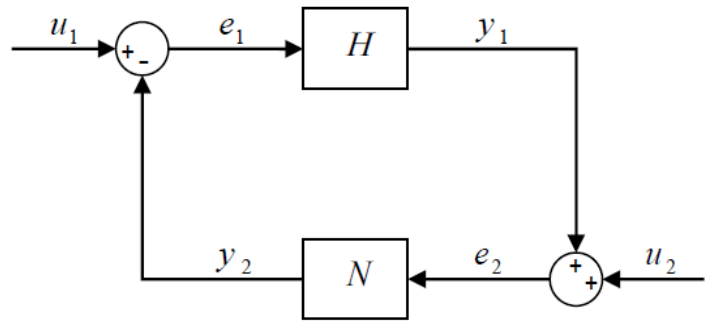

Fig. 1. A Lur'e system

Because of the generality and the wide applicability of this type of systems (for some examples, see [2], [3]), studying the stability of them has been notable and interesting to many scientists for a long time. The circle (see [4] and [5]) and Popov criteria (see [6]) can be mentioned as some of classical and important results provided yet. The approach used in the proof of Popov theorem was a functional analytic method and became a fountain of pioneering researches about the functional analysis of nonlinear systems [5], [7]-[9]. The Popov theorem is known as the first one which uses multipliers in the stability analysis of Lur'e systems. A paper of Zames and Falb published in 1968, established the multiplier method on mathematical underground [10].

Then, Cho and Narendra derived a graphical stability criterion by exploiting this fact that finding the actual multiplier is not essential [11]. One of the other input/output approaches is based on Integral Quadratic Constraint (IQC) methods that Yakubovich is known as its pioneer [12]. In 1997, Megretski and Rantzer [13] presented a stability theorem described by IQCs which is recently used in a stability theorem, referred to as generalized Zames-Falb criterion [14].

The focus of all accomplished studies of stability for Lur'e systems is on the traditional integer order systems; while, in recent decades, fractional order systems have interested many scientists in physics and engineering because of more genuine description of many actual systems made by this kind of models (for some examples, see [15-19]). One of the most important subjects in this field of research is the study of stability for fractional order systems. For the stability of linear fractional order systems, the works of Matignon [20] and Bonnet and Partington [21] can be mentioned. Momani and Hadid [22], Trasov [23] and Sabatier [24] generalized some of the stability analysis methods applied in nonlinear integer order systems to fractional order ones. Also, Ladaci and Moulay are ones who have studied of bounded input-bounded output stability for nonlinear fractional order systems [25].

In this note, input-output stability of kind $L_2$ for a fractional order system described in Lur'e structure is discussed. In this structure, the LTI subsystem considered to be of fractional order with an order between 0 and 1. To this aim, studying some of existing theorems about the stability of Lur'e systems of integer order, the conditions under which these theorems can be generalized to the corresponding fractional order systems are investigated and their applications are considered.

The paper is organized as follows. In section II, some necessary definitions and preliminaries are presented. In section III, the common classical theorems about the stability of Lur'e systems, namely circle and Popov Theorems, are investigated, and their validity for the stability of fractional order systems is studied. Then, according to the Nyquist plot of LTI systems of integer and fractional order, the application of circle criterion is compared between the fractional order Lur'e system and the integer order one. In this regard, a stability theorem is presented as well. Some multiplier theorems like Zame-Falb and generalized Zames-Falb theorem are discussed in section IV, and applying them, some classes of stable fractional order systems described in a Lur'e structure are represented. Also, another method of proof for one of the criteria used in the proof of the off-axis circle criterion is provided, so that it can be generalized to Lur'e systems of fractional order. Section V concludes the paper.

## II. Preliminaries

In this paper, $\mathbb{R}$ and $\mathbb{C}$ denotes the set of real and complex numbers, respectively. Also, let consider $L_2$ as the set of square integrable function on $\mathbb{R}$ and $L_{2e}$ as the set of functions defined on $\mathbb{R}$ which are square integrable on any compact set. $L_1$ denotes the set of absolute integrable functions on $\mathbb{R}$. In addition, we use $\hat{j}$ in this paper to denote the imaginary unit $\sqrt{-1}$. Moreover, $D[a,b]$ is a closed disk whose diameter is a line segment joins the points $(a,0)$ and $(b,0)$.

In the following, the concepts of fractional derivatives are introduced. There are different definitions for fractional derivatives; two of them are of great importance. Rienmann-Liouville derivative is defined as [26]

$$_aD_t^\alpha f(t) = \frac{d^\alpha f(t)}{d(t-a)^\alpha} = \frac{1}{\Gamma(n-\alpha)} \frac{d^n}{dt^n} \int_0^t (t-\tau)^{n-\alpha-1} f(\tau) d\tau, \tag{1}$$

where $n$ is the smallest integer order number which is greater than $\alpha$, i.e., $n-1 \leq \alpha < n$, and $\Gamma(.)$ is the Gamma function. The Caputo derivative is defined as [27]

$$_0D_t^\alpha f(t) = \begin{cases} \frac{1}{\Gamma(m-\alpha)} \frac{d^m}{dt^m} \int_0^t \frac{f^{(m)}(\tau)}{(t-\tau)^{\alpha+1-m}} d\tau, & m-1 < \alpha < m \\ \frac{d^m}{dt^m} f(t), & \alpha = m \end{cases}, \tag{2}$$

where $m$ is the first integer number which is not less than $\alpha$.

Considering all the initial conditions to be zero, the Laplace transform of Caputo fractional order derivative is [28]

$$L\left\{\frac{d^\alpha f(t)}{dt^\alpha}\right\} = s^\alpha L\{f(t)\}. \tag{3}$$

The next concept defined here is related to a fractional order transfer function which is a function of the form $H(s) = \dfrac{\sum_{i=0}^{M} b_i s^{\beta_i}}{1 + \sum_{j=1}^{N} a_j s^{\alpha_j}}$, where $0 < \alpha_1 < \alpha_2 < ... < \alpha_N$ and $0 \leq \beta_0 < \beta_1 < ... < \beta_M$. We call $H(s)$ a commensurate transfer function of order $\alpha$ if and only if it can be written it as $H(s) = F(s^\alpha)$, where $F = \dfrac{N}{D}$ that $N$ and $D$ are coprime polynomials [29]. In this paper, the LTI subsystem of the Lur'e system in the forward path is considered as a fractional order commensurate transfer function of order between 0 and 1.

The next definition is related to "Algebras" which are denoted by $A$ in this paper. Algebras are linear spaces over scalar fields $K$ with an additional property. In these spaces, a multiplication from $A \times A$ into $A$ is defined. This multiplication is associative and distributes over addition [30]. Moreover, this type of multiplication and ordinary scalar multiplication can commute. The algebra is called commutative if its special multiplication is commutative itself. If $A$ is a Banach space which is a complete normed linear space over $K$, and if $\|yz\| \leq \|y\|\|z\|$, $\forall y, z \in A$, then we call it a Banach algebra [30]. An operational example of such spaces has elements defined as [30]:

$$f(t) = \begin{cases} f_a(t) + \sum_{i=0}^{\infty} f_i \delta(t-t_i) & t \geq 0 \\ 0 & t < 0 \end{cases}, \tag{4}$$

where $f_a \in L_1[0,\infty)$; $f_i \in \mathbb{R}$ for all $i$; $\sum_{i=0}^{\infty}|f_i|<\infty$; and $t_0=0$, $t_i>0$. In this space, the product of two elements is defined to be their convolution. $\delta(t)$ which is the Dirac delta function is the unit element of this space [30].

If $G(s)$ is the Laplace transform of $g(t)$ which is an element of $A$, we say that $G$ belongs to $\hat{A}$ which is the mapping of $A$ through the Laplace transform [30].

### III. Circle and Popov Theorems for Stability Analysis of Lur'e systems

#### A. Circle Criterion

In this section, the classical form of the circle criterion for the integer order systems is stated at first, and afterwards, the proof is generalized to fractional order systems step by step. Henceforth, $g(t)$ denotes the impulse response of the LTI subsystem $H$, and $G(s)$ denotes its Laplace transform.

*Theorem III.1 (circle criterion) [30]:* Consider the system depicted in Figure 1, where $u_1,u_2,e_1,e_2,y_1,y_2:\mathbb{R}^+\to\mathbb{R}$. The subsystem $H$ is represented by

$$y_1(t)=(g*e_1)(t), \tag{5}$$

where the impulse response of $H$ in the Laplace domain is

$$G(s)=\sum_{i=1}^{k}\sum_{j=1}^{m_i}\frac{r_{ij}}{(s-p_i)^j}+G_b(s), \tag{6}$$

with $\mathrm{Re}(p_i)\geq 0$ for $i=1,...,k$, and $g_b \in L_1$ ($g_b$ is the impulse response of $G_b$). The input and output of $N$ are related with

$$y_2(t)=\Phi(e_2(t),t), \tag{7}$$

where $\Phi:\mathbb{R}\times\mathbb{R}^+\to\mathbb{R}$ belongs to the sector $\{\lambda,\gamma\}$, i.e.

$$\lambda\sigma^2 \leq \sigma\Phi(\sigma,t) \leq \gamma\sigma^2, \quad \forall t\in\mathbb{R}^+, \forall\sigma\in\mathbb{R}. \tag{8}$$

Now, let $\xi=\frac{\lambda+\gamma}{2}$, $\rho=\frac{\gamma-\lambda}{2}$, and suppose $\xi\neq 0$. With the Nyquist plot of $G$ satisfying one of the three following conditions, one can state $u_1,u_2\in L_2$ implies that $e_1,e_2,y_1,y_2\in L_2$.

a) If $0<\lambda<\gamma$, the Nyquist plot of $G$ should be out of the disk $D[-\frac{1}{\lambda},-\frac{1}{\gamma}]$ and encircle it $n_p$ times in a counterclockwise sense, where $n_p$ is the number of unstable poles of $G$.

b) If $0=\lambda<\gamma$, $G$ should not have any poles in the open right half plane and
$$\mathrm{Re}\{G(\hat{j}\omega)\}>-\frac{1}{\gamma}, \quad \forall\omega\in\mathbb{R}.$$

c) If $\lambda < 0 < \gamma$, $G$ should not have any poles in the open right half plane, and its Nyquist plot should lie entirely inside the disk $D[-\frac{1}{\lambda}, -\frac{1}{\gamma}]$.

This criterion is a direct result of the loop shifting theorem [30] and the small gain theorem [30], and it could be considered as a consequence of Theorem A.1 (see the appendix). Since the Nyquist plot of $G$ does not intersect the point $(-\frac{1}{\xi}, 0)$ and encircles it $n_p$ times in a counterclockwise direction, from Nyquist theorem, it can be stated that

$$H_\xi(s) = \frac{G(s)}{1+\xi G(s)} \in \hat{A}. \tag{9}$$

Now, letting $K = \xi$, the condition (b) of the theorem A.1 is rewritten as $\rho \sup_\omega \left|\frac{G(\hat{j}\omega)}{1+\xi G(\hat{j}\omega)}\right| < 1$. It can be easily shown that for $z \in \mathbb{C}$, the inequality $\rho|z| < |1+\xi z|$ is equivalent to one of the conditions stated in Theorem III.1.

Considering the proof of the circle criterion, in order to generalize it to fractional Lur'e systems, the two following conditions shall be met:

1) The impulse response of the subsystem $H$ belongs to $L_1$.
2) The impulse response of the subsystem $H$ belongs to $L_2$.

Condition 2 is according to Parseval's theorem [31] which is applied drawing the circle criterion.

The following two sections are devoted to studying the aforementioned conditions in fractional order systems.

### B. Belonging of the Impulse Response of Fractional Order Systems to $L_1$ Space

Belonging of the impulse response of fractional order systems to $L_1$ space is discussed in this section.

*Definition III.1 [20]:* A linear system with the impulse response $h$ is $L_\infty$-BIBO stable iff $\forall u \in L_\infty \Rightarrow y = h * u \in L_\infty$. This definition is equivalent to $h \in L_1$.

*Theorem III.2 [32]:* A commensurate transfer function of order $\alpha$ which is in the form of $H(s^\alpha) = \frac{N(s^\alpha)}{D(s^\alpha)}$ where $N$ and $D$ are coprime polynomials is $L_\infty$-BIBO stable iff $0 < \alpha < 2$ and for every $s \in \mathbb{C}$ where $D(s) = 0$ for $|\arg(s)| > \alpha \frac{\pi}{2}$.

As a result, if a fractional order system satisfies the conditions of Theorem III.2, then its impulse response is in $L_1$. As some examples, we can mention systems with transfer functions of the form

$$G(s^\alpha) = \sum_{i=1}^{k}\sum_{j=1}^{m_i} G_{ij}(s^\alpha), \quad G_{ij}(s^\alpha) = \frac{r_{ij}}{(s^\alpha + p_i)^j}, \tag{10}$$

where $p_i > 0$. From Theorem III.2, these systems have impulse responses with finite $L_1$-norms.

## C. Belonging of the Impulse Response of Fractional Order Systems to $L_2$ Space

The following theorem discusses the conditions under which the boundedness of $L_2$-norm of impulse response of a fractional order system is guaranteed.

*Theorem III.3 [29]:* Let $F(s^\alpha)$ be a commensurate transfer function of order $\alpha$ which is $L_\infty$-BIBO stable. Additionally, assume that the numerator and denominator of $F(s^\alpha)$ are coprime. Let $\beta_M$ and $\alpha_N$ be the numerator and denominator degree of $F(s^\alpha)$, respectively. The $L_2$-norm of impulse response of $F(s^\alpha)$ is bounded if $\alpha_N - \beta_M > \frac{1}{2}$.

Now, we can state a corollary which states some conditions to generalize the circle criterion to the Lur'e systems of fractional.

*Corollary III.1:* Consider the Lur'e system shown in Figure 1. Suppose that the LTI subsystem $H$ is of fractional order and has a commensurate transfer function of order $\alpha$ which is in the form of $G(s^\alpha) = \frac{N(s^\alpha)}{D(s^\alpha)}$ where $N$ and $D$ are coprime polynomials, and $0 < \alpha \le 1$. Let $\beta_M$ and $\alpha_N$ be the numerator and denominator degree of $G(s^\alpha)$, respectively. If impulse response $H$ satisfies the necessary conditions of being in $L_1$ and $L_2$, i.e. 1) for every $s \in \mathbb{C}$ where $D(s) = 0$, $|\arg(s)| > \alpha\frac{\pi}{2}$, and 2) $\alpha_N - \beta_M > \frac{1}{2}$, then the circle criterion can be generalized to the related fractional order Lur'e system.

## D. The Application of the Circle Criterion in the Integer and Fractional Order Lur'e Systems

One of significant advantages of the circle criterion is its easily utility by drawing the Nyquist plot of the linear subsystem. Therefore, in this section, we intend to compare the Nyquist plots of integer order system $G(s)$ and its corresponding fractional order system $G(s^\alpha)$ with $0 < \alpha \le 1$.

The next definition and lemma define a strictly positive real transfer function and propose an equivalent description for it.

*Definition III.2 [33]:* A real rational transfer function $H(s)$ is called positive real (PR) if $\text{Re}\{H(s)\} \ge 0$ for all $s$ with a nonnegative real part. The transfer function $H(s)$ is called strictly positive real (SPR) if $H(s - \varepsilon)$ is positive real for some $\varepsilon > 0$.

*Lemma III.3 [33]:* The proper real transfer function $H(s)$ is positive real iff the following conditions are satisfied.

1. No poles of $H(s)$ have positive real parts.
2. If some poles of $H(s)$ are on the imaginary axis, then they should be simple poles with positive residues.
3. $\text{Re}\{H(\hat{j}\omega)\} \geq 0$.

In addition, if all the poles of $H(s)$ are in the left half plane, and $\text{Re}\{H(\hat{j}\omega)\} > 0$, $H(s)$ is called strictly positive real.

The next theorem compares the Nyquist plots of the system $G(s)$ and its corresponding fractional order system $G(s^\alpha)$ with $0 < \alpha \leq 1$.

*Theorem III.4:* Let $G(s)$ be an integer order proper transfer function which satisfies the following conditions.

1. $G(s)$ has no poles in the right half plane.
2. If some poles of $G(s)$ are on the imaginary axis, then they should be simple poles with positive residues.
3. The Nyquist plot of $G(s)$ lies to the right of the vertical line $\text{Re}\{s\} = a$ ; i.e. $\text{Re}\{G(\hat{j}\omega)\} \geq a, \quad \forall \omega \in \mathbb{R}$.

Now, if the Nyquist plot of $G(s)$ lies to the right of the vertical line $\text{Re}\{s\} = a$, then it can be deduced that the Nyquist plot of $G(s^\alpha)$ with $0 < \alpha \leq 1$ lies to the right of the same line as well. Also, if $G(s)$ has no poles on the imaginary axis and its Nyquist plot lies to the left of vertical line $\text{Re}(s) = b$, then the Nyquist plot of $G(s^\alpha)$ with $0 < \alpha \leq 1$ lies to the left of the same line, too.

*Proof:* According to the assumptions, it can be stated that $\text{Re}\{G(\hat{j}\omega)\} \geq a, \quad \forall \omega \in \mathbb{R}$. From Lemma III.3, it can be shown that $H(s) = G(s) - a$ is positive real; i.e.

$$\text{Re}\{H(s)\} \geq 0, \quad \forall s, \quad \text{Re}\{s\} \geq 0. \tag{11}$$

Now, consider the mapping $s^\alpha$ with $\alpha = \frac{u}{v} \in \mathbb{Q}$. This mapping is multivalued so that for $s = re^{j\theta}$ and $-\pi < \theta \leq \pi$, we have $s^{\frac{u}{v}} = r^{\frac{u}{v}} e^{\hat{j}(\frac{u\theta + 2k\pi}{v})}$, $k = 0, 1, \ldots, v - 1$. The principal value of $s^\alpha$ is obtained by $k = 0$ and has a physical interpretation [34]. Then, it can be stated that

$$\text{Re}\{\hat{j}^\alpha \omega^\alpha\} = \omega^\alpha \cos\left(\frac{\alpha\pi}{2}\right) \geq 0, \quad \forall \omega \in \mathbb{R}. \tag{12}$$

Thus, by (11), we obtain

$$\text{Re}\{H(\hat{j}^\alpha \omega^\alpha)\} \geq 0, \quad \forall \omega \in \mathbb{R}. \tag{13}$$

Since $H(\hat{j}^\alpha \omega^\alpha) = G(\hat{j}^\alpha \omega^\alpha) - a$, it can be deduced from (13) that

$$\text{Re}\{G(\hat{j}^\alpha \omega^\alpha)\} \geq a, \quad \forall \omega \in \mathbb{R}. \tag{14}$$

Therefore, the Nyquist plot of $G(s^\alpha)$ lies to the right of the vertical line $\text{Re}\{s\}=a$ as well. The proof of the second case in which the Nyquist plot of $G(s)$ and $G(s^\alpha)$ lie to the left of $\text{Re}(s)=b$, is rather similar and is omitted.

*Corollary III.2*: Suppose that the fractional-order Lur'e system is defined as the negative feedback interconnection of a nonlinear subsystem $N$ and a LTI subsystem with the transfer function $G(s^\alpha)$, where $0 < \alpha \leq 1$, which satisfies the conditions of the circle criterion. Also, assume that the corresponding integer order transfer function $G(s)$ has no poles in the closed right half plane and its Nyquist plot lies to the right of the vertical line $\text{Re}\{s\} = -\dfrac{1}{\beta}$. Thus, according to the circle criterion and Theorem III.4, it can be derived that the fractional order Lur'e system is $L_2$-stable if $N$ belongs to the sector $\{0,\beta\}$.

As a numerical example of theorem III.4, consider the Nyquist plot of the integer-order system $G(s) = \dfrac{6}{(s+1)(s+2)(s+3)}$ and its corresponding fractional order system $G(s^{0.7}) = \dfrac{6}{(s^{0.7}+1)(s^{0.7}+2)(s^{0.7}+3)}$. The Nyquist plots of the two systems are depicted in Figure 2. The Nyquist plot of $G(s)$ lies to the right of $\text{Re}\{s\}=-0.2147$, and the Lur'e system defined as the negative feedback interconnection of $G(s)$ and the nonlinear subsystem which belongs to the sector $\{0, 4.6577\}$ is $L_2$-stable.

In addition, the Nyquist plot of $G(s^{0.7})$ lies to the right of $\text{Re}\{s\}=-0.03745$, and $L_2$-stability is guaranteed for a larger sector $\{0, 26.7023\}$.

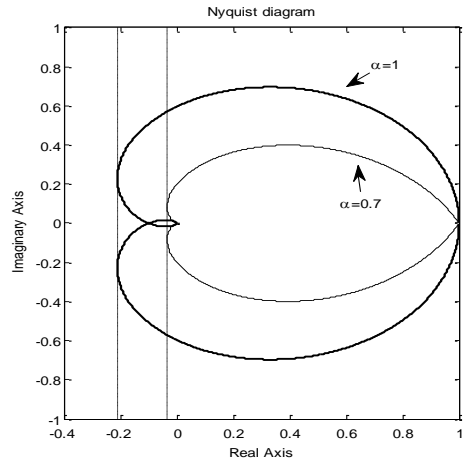

Fig. 2. Nyquist plots of the transfer functions $G(s) = \dfrac{6}{(s+1)(s+2)(s+3)}$ and $G(s) = \dfrac{6}{(s^\alpha+1)(s^\alpha+2)(s^\alpha+3)}$

In order to test the validity of the obtained results, we have simulated two fractional and integer order Lur'e systems whose LTI subsystems are described with the transfer functions $G(s^\alpha) = \dfrac{6}{(s^\alpha+1)(s^\alpha+2)(s^\alpha+3)}$, where

$\alpha = 1$, 0.7. Also, the nonlinear subsystem is a saturation function which is defined as:

$$(Ne_2)(t) = \Phi[e_2(t)] = \begin{cases} 1, & e_2(t) \geq 1 \\ e_2(t), & |e_2(t)| < 1 \\ -1, & e_2(t) \leq -1 \end{cases}, \tag{15}$$

Note that the nonlinearity belongs to the permissible sectors for both systems. In these simulations, an Adams-type predictor-corrector method [35] for the numerical solution of fractional differential equations is applied. The input $u_1$ is a square pulse which is defined as:

$$u_1 = \begin{cases} 5 & 0 \leq t \leq 50 \\ 0 & oth. \end{cases}. \tag{16}$$

The time responses of the output $y_1$ are illustrated in Figure 3.

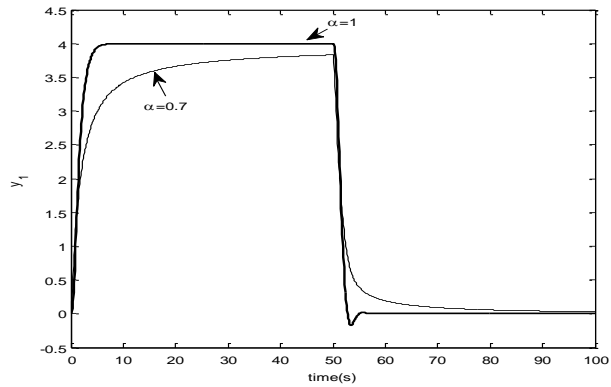

Fig. 3. Pulse responses of the Lur'e systems which are interconnections of the saturation nonlinearity with the description (15) in the feedback path and the LTI systems $G(s^\alpha) = \dfrac{6}{(s^\alpha + 1)(s^\alpha + 2)(s^\alpha + 3)}$, where $\alpha = 1$, 0.7 in the forward path.

As another numerical example, we consider the same LTI subsystems as the previous example, but we assume that the nonlinearity belong only to the permissible sector of the fractional order Lur'e system with $\alpha = 0.7$. It is defined as:

$$(Ne_2)(t) = \Phi[e_2(t)] = \begin{cases} 10, & e_2(t) \geq 1 \\ 10e_2(t), & |e_2(t)| < 1 \\ -10, & e_2(t) \leq -1 \end{cases} \tag{17}$$

The input $u_1$ is a square pulse whose value is 20 in the time interval $\left[0^{\sec}, 50^{\sec}\right]$, and otherwise it is equal to zero. The responses are depicted in Figure 4. It is obvious that the output of the integer order Lur'e system is not in $L_2$; thus, this system cannot be $L_2$ stable, while the corresponding fractional order one can be stable.

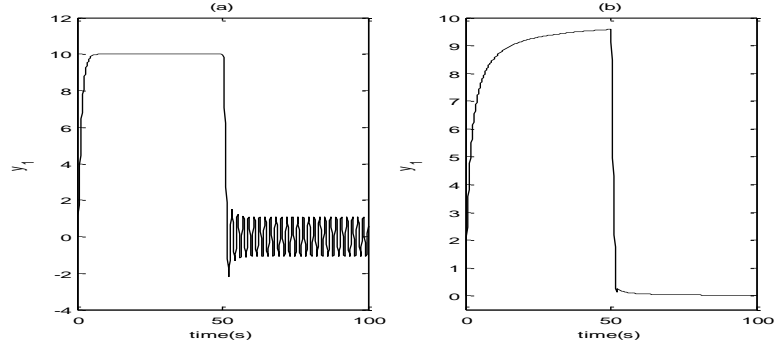

Fig. 4. Pulse responses of the Lur'e systems which are interconnections of the saturation nonlinearity with the description (17) in the feedback path and the LTI systems $G(s^\alpha) = \dfrac{6}{(s^\alpha+1)(s^\alpha+2)(s^\alpha+3)}$ in the forward path, where (a) $\alpha = 1$, and (b) $\alpha = 0.7$.

### E. Popov Theorem

In this section, first, the Popov criterion is proposed, and then the conditions under which this criterion can be generalized to fractional-order Lur'e systems are discussed.

*Theorem III.5 [30]:* Consider a single input-single output system which is described by

$$\begin{aligned} e_1 &= u_1 - g * e_2 \\ e_2 &= \Phi(e_1) \end{aligned}, \quad (18)$$

where $g \in L_1[0,\infty)$, $\dot{g} \in A$, $\Phi: \mathbb{R} \to \mathbb{R}$ is continuous and belongs to the sector $\{0, k\}$. Assume that for every $u_1 \in L_2$, we have that $e_1, e_2 \in L_{2e}$. If for some $q \geq 0$, there is a $\varepsilon > 0$ so that

$$\text{Re}\left\{(1+q\hat{j}\omega)G(\hat{j}\omega)+\frac{1}{k}\right\} \geq \varepsilon > 0, \quad \forall \omega \geq 0, \quad (19)$$

then, for every $u_1, \dot{u}_1 \in L_2$, $e_1, \dot{e}_1, e_2 \in L_2$. Also, $e_1$ belongs to $L_\infty$ space and decays toward zero as $t \to \infty$.

This theorem can be proved using the small gain theorem [30], the passivity theorem [30], and the loop transformation theorem [30]. For a complete proof, refer to [30].

From the proof of Popov's theorem, we can deduce that if the LTI subsystem $G$ is of fractional order, its impulse response $g$ should be in $L_1$ and $L_2$, and $\dot{g} \in A$.

*Corollary III.3:* Consider the Lur'e system depicted in Figure 1. Assume that the transfer function of the LTI subsystem $H$ is of fractional-order and described by

$$G(s^\alpha) = \dfrac{\sum_{i=0}^{m} b_i s^{i\alpha}}{\sum_{j=0}^{n} a_j s^{j\alpha}}, \quad b_m, a_n \neq 0. \quad (20)$$

In order to generalize the Popov theorem to this kind of Lur'e systems, the mentioned necessary conditions for belonging $G$ to $L_1$ and $L_2$ should be satisfied. Also, we should have $\dot{g} \in A$. If $g(0) = 0$, then the Laplace transform of $\dot{g}$ is $sG(s^\alpha) = \dfrac{\sum_{i=0}^{m} b_i s^{i\alpha+1}}{\sum_{j=0}^{n} a_j s^{j\alpha}}$. If $n\alpha \geq m\alpha + 1$, then $\dot{g} \in A$. In other words, the relative degree of the LTI subsystem should not be less than the inverse of $\alpha$.

## IV. Stability Analysis of Lur'e Systems Using Multipliers and Integral Quadratic Constraints

In this section, some theorems deal with stability of Lur'e systems using multipliers, are stated and the necessary conditions to apply them in the fractional order Lur'e systems are discussed. Also, some classes of fractional order Lur'e systems whose stability can be proved using the aforementioned theorems are introduced.

Perhaps, one of the first stability theorems using multipliers is the Popov's theorem. In fact, the main idea in multiplier-type stability theorems is satisfying the passivity or positivity conditions by multiplying the LTI operator in an appropriate chosen multiplier. In the following, the Zames-Falb multiplier theorem is stated.

*Theorem IV.1 (Zames-Falb Theorem) [10]:* Consider the Lur'e system illustrated in Figure 1 with the following equation holds:

$$y_1(t) = He_1(t) = \int_0^\infty g(\tau) e_1(t-\tau) d\tau, \tag{21}$$

where $g \in L_1[0, \infty)$. Let $G(\hat{j}\omega)$ denote the frequency response of the LTI subsystem $H$. $N$ is the nonlinear subsystem such that $y_2(t) = Ne_2(t) = \Phi(e_2(t))$, where the following assumptions should be met:

1) $\Phi(0) = 0$.
2) There exists a constant $C > 0$, where $|\Phi(r)| \leq C|r|$, for all real $r$.
3) $\Phi$ is monotone non-decreasing; i.e. $(r-s)[\Phi(r) - \Phi(s)] \geq 0$, for all real $r$ and $s$.

An additional trivial assumption is that for all $u_1, u_2 \in L_2$, $e_1, e_2 \in L_{2e}$.

Now, if there exists a LTI operator $z$ with the frequency response $Z(\hat{j}\omega)$ such that

1) $\int_{-\infty}^{\infty} |z(t)| dt < 1,$ \qquad (22)

2) $\text{Re}\left\{\left[1 - Z(\hat{j}\omega)\right] G(\hat{j}\omega)\right\} \geq \varepsilon > 0, \quad \forall \omega \in \mathbb{R},$ \qquad (23)

and if either $z(\cdot) \geq 0$, or $\Phi(\cdot)$ is odd, then $e_1$ and $e_2$ are in $L_2(-\infty, \infty)$. In addition, if $g \in L_2[0, \infty)$, then $y_1(t)$ decays towards zero as $t \to \infty$.

In the Zames-Falb theorem, a general structure to use multipliers in the stability analysis of a system is proposed. This theorem is a result of the passivity theorem and a multiplier theorem which uses noncausal multipliers to prove the stability of a feedback system. The conditions of this theorem are derived from some positivity lemmas and a factorization theorem which deals with factorizing an operator into the product of causal and noncausal operators. Since this theorem is an input-output stability theorem, and uses functional analytic approach to prove stability, it can be applied in fractional order systems as well. Indeed, in functional analytic method, the systems are considered as operators, and it helps to its utility in system of fractional order. A precisely surveillance of all details of the proof, which applies Parseval's theorem and the properties of the Banach algebra of convolution operators [30], results in the following corollary.

*Corollary IV.1:* if $L_1$ and $L_2$-norm of the impulse response of the LTI subsystem $H$ are bounded, applying the Zames-Falb theorem in the fractional order Lur'e systems is possible.

**A. Some Classes of Stable Fractional-Order Lur'e Systems**

In this part, applying the Zames-Falb theorem, we want to introduce some classes of stable fractional order Lur'e systems. To this aim, let us review the important properties of $Z$ in the Zames-Falb theorem at first. Indeed, $\|Z\|_1$ should be less than one, and if $\Phi(\cdot)$ is not an odd function, $Z$ should have a nonnegative impulse response. The following lemma, relates nonnegativity and $L_1$-norm of the fractional order $Z(s^\alpha)$ to the integer order $Z(s)$.

*Lemma IV.2:* If the impulse response of a BIBO-stable integer order transfer function $Z(s)$ is nonnegative, then the fractional order transfer function $Z(s^\alpha)$, where $0<\alpha<1$, has a nonnegative impulse response as well, and its $L_1$-norm is equal to $Z(0)$.

*Proof:* If the impulse response of $Z(s)$ is nonnegative, its step response is monotone non-decreasing, and from Theorem 3 of [36], the fractional order $Z(s^\alpha)$ has a monotone non-decreasing step response too. Hence, its impulse response will be nonnegative. For a transfer function $Z$ with a nonnegative impulse response, we have

$$\|Z\|_1 = \int_0^\infty |z(t)|dt = \int_0^\infty z(t)dt \tag{24}$$

Thus, $\|Z\|_1$ is the final value of step response of $Z$, and can be computed using final value theorem as $\|Z\|_1 = Z(s)|_{s=0}$ [37]. ∎

The next theorem which is deduced from Lemma IV.2, Theorem III.4, and Zames-Falb theorem proves the stability of a fractional order Lur'e system, if the corresponding integer order Lur'e system is stable.

*Theorem IV.3:* If for an integer order Lur'e system which is depicted in Figure 1, the stability is proved by the Zames-Falb theorem, and the integer order $Z(s)$ which satisfies (22) and (23), is BIBO stable, and it has a nonnegative impulse response with the $L_1$-norm less than 1, then the corresponding fractional order Lur'e system with the LTI subsystem $G(s^\alpha)$, where $0<\alpha<1$, is stable as well.

*Proof:* Let us define $F(s)=[1-Z(s)]G(s)$. If (23) is satisfied, from the Theorem III.4, we have that $\text{Re}\{F(\hat{j}^\alpha \omega^\alpha)\} \geq \varepsilon > 0$. Additionally, using Lemma IV.2, the impulse response of $Z(s^\alpha)$ is nonnegative, with the $L_1$-norm equal to $Z(0)$ which is less than 1. Thus, from the Zames-Falb theorem, the stability of the corresponding fractional order Lur'e system is proved. ∎

Now, using the transfer functions with nonnegative impulse responses listed in [36], we can introduce some classes of transfer functions $G(s)$ which satisfy (23) and result in stable fractional order Lur'e systems.

Note that by adding or convolving the nonnegative functions, more fractional order transfer functions with nonnegative impulse responses can be generated.

1) $$G(s) = k + \sum_{i=1}^{n} \frac{k_i}{s^{\alpha_i} + b_i} \tag{25}$$

where $k, k_i > 0$, $b_i > 0$ and $0 < \alpha_i \leq 1$. It is obvious that $\text{Re}\{G(\hat{j}\omega)\} \geq \varepsilon > 0$, and by choosing $Z(s) \equiv 0$, (23) is satisfied.

2) $$G(s) = \left(\frac{s^\alpha + b}{s^\alpha + a}\right)\left(k + \sum_{i=1}^{n} \frac{k_i}{s^{\alpha_i} + b_i}\right) \tag{26}$$

where $k, k_i > 0$, $b > a > 0$, $b_i > 0$ and $0 < \alpha, \alpha_i \leq 1$. Choosing $Z(s) = \frac{b-a}{s^\alpha + b}$, (23) is satisfied. Also, we have $\|Z\|_1 = Z(0) = \frac{b-a}{a} < 1$.

3) $$G(s) = \frac{(s^\alpha + a)(s^\alpha + c)}{\left[s^{2\alpha} + (a+c-k_0)s^\alpha + ac - k_0 b\right]}\left(k + \sum_{i=1}^{n} \frac{k_i}{s^{\alpha_i} + b_i}\right) \tag{27}$$

where $k, k_0, k_i, c > 0$, $b > a > 0$, $kb < ac$, $a + c > k$, $b_i > 0$ and $0 < \alpha, \alpha_i \leq 1$.

In this case, let $Z(s)$ be defined as $Z(s) = \frac{k(s^\alpha + b)}{(s^\alpha + a)(s^\alpha + c)}$ whose $L_1$-norm is equal to $\frac{kb}{ac}$ which is less than 1. $G(s)$ is obtained by multiplying $\frac{1}{1-Z}$ and a transfer function with a positive real part on the imaginary axis, which satisfies (23).

Now, Consider Lur'e systems as depicted in Figure 1, with the nonlinear subsystem (15) and the LTI subsystems described as:

$$G(s) = 10^{-6} + \frac{1}{s^\alpha + 1} + \frac{2}{s^\alpha + 2}, \tag{28}$$

$$G(s) = \frac{s^\alpha + 4}{s^\alpha + 3}\left(10^{-6} + \frac{1}{s^\alpha + 1} + \frac{2}{s^\alpha + 2}\right), \tag{29}$$

$$G(s) = \frac{(s^\alpha + 3)(s^\alpha + 5)}{s^{2\alpha} + 8s^\alpha + 7}\left(10^{-6} + \frac{1}{s^\alpha + 1} + \frac{2}{s^\alpha + 2}\right), \tag{30}$$

which are respectively the numerical examples of case (1), (2), and (3) introduced above, and $\alpha = 1, 0.7$. The pulse responses of the Lur'e systems to the input (16) are shown in Figure 5. It is obvious that the outputs of the systems decay toward zero as the Zames-Falb theorem predicts.

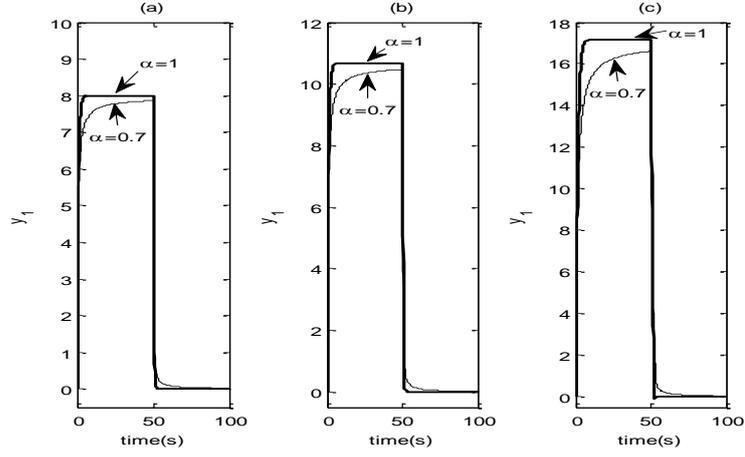

Fig. 5. Pulse responses of the Lur'e systems which are interconnections of the saturation nonlinearity with the description (15) in the feedback path and the LTI subsystems in the forward path described by (a) the relation (29), (b) the relation (30), and (c) the relation (31), where $\alpha = 1$, $0.7$.

### B. The Generalized Zames-Falb Multiplier Theorem

The Generalized Zames-Falb theorem is based on integral quadratic constraints (IQC's) which is a unifying framework for robust control problems. This criterion provides robust stability conditions which consider quasi-monotone and odd behaviors. In this section, the theorem is stated, and then a class of stable fractional order Lur'e systems, utilizing it, is introduced.

*Theorem IV.4 (Generalized Zames-Falb) [14]:* Consider the Lur'e system illustrated in Figure 1. Assume that $G(s)$ is stable, and $y_2(t) = n(e_2(t))$, such that $n$ is quasi-monotonic-and-odd, i.e. there exists a function $\delta : \mathbb{R} \to \mathbb{R}$ such that $n(y) = \bar{n}(y)[1+\delta(y)]$, $\bar{n}:\mathbb{R} \to \mathbb{R}$ is monotone non-decreasing and odd, and $|\delta(y)| \leq D < 1$. Moreover, the feedback interconnection of $G$ and $\tau N$ is well-posed for every $\tau \in [0,1]$, as defined in [14]. In addition, , $N$ is bounded, $yn(y) \geq 0$, and there exists $\varepsilon > 0$ and a function $z \in L_1$ such that

$$\|z\|_1 \leq \left(\frac{1-D}{1+D}\right)^2, \tag{31}$$

$$\mathrm{Re}\left\{G(\hat{j}\omega)\left[1-Z^*(\hat{j}\omega)\right]\right\} \geq \varepsilon G^*(\hat{j}\omega)G(\hat{j}\omega), \tag{32}$$

where $Z(j\omega)$ is the frequency response of $z(t)$. Then, the Lur'e system is $L_2$-stable.

The proof of this theorem is based on a stability theorem by Megretski and Rantzar [14] who have used IQC's in the stability analysis. In this proof, the Parseval's theorem and the properties of norms are utilized. Thus, the following result is obtained.

*Corollary IV.5:* The boundedness of $L_1$ and $L_2$-norm of the impulse response for the LTI subsystem $G(s^\alpha)$ in a fractional order Lur'e system is the necessary condition to apply the generalized Zames-Falb theorem in such systems.

Using the generalized Zames-Falb theorem, a class of stable fractional order Lur'e systems can be introduced. Since $G$ has a bounded $L_1$-norm, and none of its zeros are on the imaginary axis, we have that

$$\text{Re}\left\{G(\hat{j}\omega)\left[1-Z^*(\hat{j}\omega)\right]\right\} = |G(\hat{j}\omega)|^2 \text{Re}\left\{\frac{1-Z^*(\hat{j}\omega)}{G^*(\hat{j}\omega)}\right\}.$$ Thus, the relation (32) is equivalent to

$$0 < \text{Re}\left\{\frac{G(\hat{j}\omega)}{1-Z(\hat{j}\omega)}\right\} \leq \text{E}. \tag{33}$$

where $\text{E} \in \mathbb{R}$ is a bounded number. With respect to the discussions in section IV.A, a function $Z$ can be defined as

$$Z = K \frac{\prod_{i=1}^{m}(s^\alpha + b_i)}{\prod_{j=1}^{n}(s^\alpha + a_j)}, \tag{34}$$

where $n \geq m$, $0 < \alpha \leq 1$, and for every $b_i$, there is a $a_j$ such that $b_i > a_j$. Moreover, we should have that

$\|Z\|_1 = K \dfrac{\prod_{i=1}^{m} b_i}{\prod_{j=1}^{n} a_j} < 1$. Then, $G$ can be defined as a multiplication of $1-Z$ and a transfer function with a positive bounded real part on the imaginary axis. A class of fractional order LTI subsystems which result in stable Lur'e systems and satisfy the conditions of generalized Zames-Falb theorem, can be defined as

$$G(s) = \left(1 - K \frac{\prod_{i=1}^{m}(s^\alpha + b_i)}{\prod_{j=1}^{n}(s^\alpha + a_j)}\right)\left(\sum_{l=1}^{n} \frac{k_l}{s^{\alpha_l} + b_l}\right), \tag{35}$$

where $k_l > 0$, $b_l > 0$, $0 < \alpha_l \leq 1$, and for every $b_i$, there is a $a_j$ such that $b_i > a_j$.

As a numerical example, consider the Lur'e system depicted in Figure 1 with the nonlinear subsystem (15) and the LTI subsystem described as:

$$G(s) = \left(1 - 0.5\frac{s^\alpha + 4}{s^\alpha + 3}\frac{s^\alpha + 6}{s^\alpha + 5}\right)\left(\frac{1}{s^\alpha + 1} + \frac{2}{s^\alpha + 2}\right), \tag{36}$$

where $\alpha = 1$ and $0.7$. The pulse responses of the Lur'e systems to the input (16) are depicted in Figure 6.

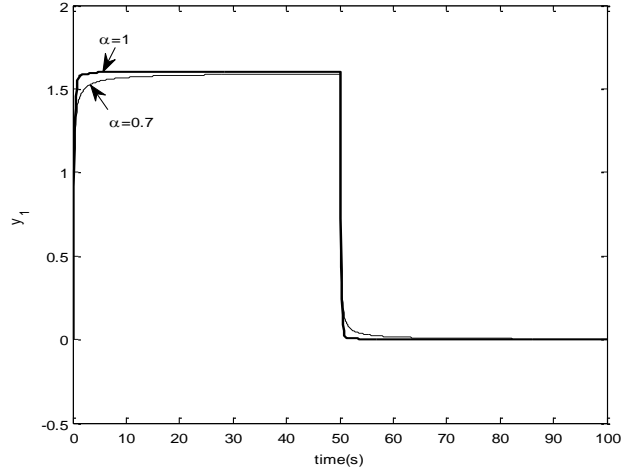

Fig. 6. Pulse responses of the Lur'e systems which are interconnections of the saturation nonlinearity with the description (15) in the feedback path and the LTI subsystems in the forward path described by the relation (37), where $\alpha = 1$, $0.7$.

### C. A New Verification Method of Skeleton Theorem

One of the important properties of the multiplier-type stability techniques is that ensuring the existence of an appropriate multiplier is sufficient, and finding it is not necessary. The off-axis circle criterion [11] is a method of stability analysis which guarantees the existence of a multiplier without finding it via some geometric conditions. A criterion named skeleton theorem [8] is utilized in derivation of the off-axis circle criterion. However, its proof proposed by Zames is not extendable to the fractional order Lur'e systems. In this section, we want to develop another proof for this theorem, so that the off-axis circle criterion can be generalized to Lur'e systems of fractional order. In this proof, we use the Zames-Falb and the Generalized Zames-Falb theorem mentioned before. First, let us define RL class of functions as

$$M(s) = \prod_{i=1}^{N} \frac{s + \alpha_i}{s + \beta_i}, \tag{37}$$

where $0 < \alpha_1 < \beta_1 < \alpha_2 < \beta_2 < \ldots < \alpha_N < \beta_N$. In addition, RC is a class of function $M(s)$ for which $M^{-1}(s)$ is in RL.

*Theorem IV.6 (Skeleton Theorem) [8]:* Consider the Lur'e system depicted in Figure 1, and assume that $u_2 \equiv 0$. In addition, $g(t)$ is in $L_1$ and $L_2$, $N$ is nonlinear and time-invariant, and $N(0) = 0$. For every $y_1, y_2 \in \mathbb{R}$ such that $y_1 \neq y_2$, we have that

$$K_1 \leq \frac{N(y_1) - N(y_2)}{y_1 - y_2} \leq K_2. \tag{38}$$

If there exists a multiplier $M(s)$ in RL or RC such that

$$\operatorname{Re}\left\{ \frac{K_2 G(\hat{j}\omega) + 1}{K_1 G(\hat{j}\omega) + 1} M(\hat{j}\omega) \right\} \geq \varepsilon > 0, \quad \omega \in \mathbb{R}, \tag{39}$$

then $u_1(t) \in L_2$ implies $y_1(t) \in L_2$, and $\lim_{t \to \infty} y_1(t) = 0$.

To derive this theorem in the case that $M(s)$ is in RC, Zames assumed that $H$ can be factorized as $H = KH_1$ such that $K$ is in RC. Then, he synthesized $N$ and $K$ in the closed loop system and proved that if $N$ is extended incrementally positive, as defined in [8], $NK$ is extended incrementally positive as well. This factorization approach is suitable for the integer order Lur'e systems, but it cannot be easily extended to the fractional order ones. In the following, another verification method of the skeleton theorem is proposed, which can be generalized to the fractional order systems.

*Proof of Theorem IV.6:* Applying the loop transformation theorem [30], we prove the skeleton criterion for the case that the nonlinear subsystem $N$ is monotone non-decreasing, and $N(0) = 0$. In this case, the relation (39) is rewritten as

$$\operatorname{Re}\{G(\hat{j}\omega)M(\hat{j}\omega)\} \geq \varepsilon > 0. \tag{40}$$

Now, we want to show that the skeleton is a special form of the Zames-Falb theorem if $M(s)$ is RL. To this aim, we prove that $M(s)$ can be considered as $1 - Z(s)$, where $\|z\|_1 < 1$. A RL function $M(s)$, as defined in (38), can be rewritten as

$$M(s) = 1 + \sum_{i=1}^{N} \frac{k_i}{s + \beta_i} \tag{41}$$

where $k_i = (\alpha_i - \beta_i) \prod_{\substack{j=1 \\ j \neq i}}^{N} \frac{\alpha_j - \beta_i}{\beta_j - \beta_i}$, $i = 1, \ldots, N$. Since $\alpha_i < \beta_i$, $\alpha_i - \beta_i$ is negative. Furthermore, for $j \neq i$, since $\alpha_j$ and $\beta_j$ are both less than or greater than $\beta_i$, then $\frac{\alpha_j - \beta_i}{\beta_j - \beta_i}$ is positive. Thus, $k_i$'s are negative and $Z(s)$ can be defined as $Z(s) = \sum_{i=1}^{N} \frac{-k_i}{s + \beta_i}$ which has a nonnegative impulse response. Since $z(t) \geq 0$, From Lemma IV.2, $0 < \|z\|_1 = Z(0) = 1 - \prod_{i=1}^{N} \frac{\alpha_i}{\beta_i} < 1$. Hence, (40) is equivalent to (23), and the $L_2$-stability of the Lur'e system is resulted by Zames-Falb theorem.

Now, consider the case that $M(s)$ is a RC function. In this case, $M^{-1}(s) = 1 - Z(s)$ is a RL function, where $\|z\|_1 < 1$. Thus, if there exists a RC $M(s)$ which satisfies (40), (33) holds as well, and the stability is proved by the generalized Zames-Falb theorem. ∎

## V. Conclusion

In this paper, the $L_2$ input-output stability of a fractional order system which is described in a Lur'e structure has been investigated. We have concluded that the impulse response of the LTI subsystem of a fractional order Lur'e system with an order between 0 and 1 should have bounded $L_1$ and $L_2$ norms so that the circle and Popov criteria can be applied in the stability analysis. Furthermore, comparing the Nyquist plots of LTI systems of integer order and fractional order, it has been shown that if a Lur'e system of integer order is stable, the corresponding Lur'e system of fractional order with an order in (0,1] can be stable as well, and it may be stable even for a greater sector that the nonlinear subsystem belongs to. In part IV, some of the multiplier-type stability theorems have been studied and some classes of stable fractional order Lur'e systems have been introduced. Finally, in order to generalize the off-axis circle criterion to the Lur'e systems of fractional order, another approach is presented to prove the skeleton theorem used in its derivation.

**Appendix**

*Theorem A.1 [30]:* Consider the system shown in Figure 1. The subsystem $H$ is defined by $(He_1)(t) = \int_0^t g(t-\tau)e_1(\tau)d\tau$. $G(s)$ stands for the Laplace transform of $g$, and $G(s) = N(s)D^{-1}(s)$ such that $N$ and $D$ are coprime polynomials. The subsystem $N$ is represented by $(Ne_2)(t) = \Phi(e_2(t),t)$. If there exist real constants $K$ and $\gamma$ so that $(Ne_2)(t) = \Phi(e_2(t),t)$, and if

a) $\inf_{Re(s) \geq 0} |1 + KG(s)| > 0$,

b) $\gamma \sup_{\omega \in \mathbb{R}} |H_K(\hat{j}\omega)| < 1$, where $H_K(s) = G(s)[1 + KG(s)]^{-1}$

Then, $u_1, u_2 \in L_2 \Rightarrow e_1, e_2, y_1, y_2 \in L_2$.

Loop shifting theorem [30] can be utilized in the proof of this theorem. The condition (a) is used to show that $H_k \in \hat{A}$. The result is finally obtained, from small gain theorem.